# A Deep Reinforcement Learning Based Approach for Optimal Active Power Dispatch

Jiajun Duan[1], Haifeng Li[2], Xiaohu Zhang[1], Ruisheng Diao[1], Bei Zhang[1],
Di Shi[1*], Xiao Lu[2], Zhiwei Wang[1], Siqi Wang[1]

[1]Global Energy Interconnection Research Institute North America (GEIRINA), San Jose, USA
[2]State Grid Jiangsu Electric Power Company, Nanjing, CN
di.shi@geirina.net

*Abstract*—The stochastic and dynamic nature of renewable energy sources and power electronic devices are creating unique challenges for modern power systems. One such challenge is that the conventional mathematical systems models-based optimal active power dispatch (OAPD) method is limited in its ability to handle uncertainties caused by renewables and other system contingencies. In this paper, a deep reinforcement learning-based (DRL) method is presented to provide a near-optimal solution to the OAPD problem without system modeling. The DRL agent undergoes offline training, based on which, it is able to obtain the OAPD points under unseen scenarios, e.g., different load patterns. The DRL-based OAPD method is tested on the IEEE 14-bus system, thereby validating its feasibility to solve the OAPD problem. Its utility is further confirmed in that it can be leveraged as a key component for solving future model-free AC-OPF problems.

*Index Terms*—Artificial intelligence, optimal active power dispatch, deep reinforcement learning, neural network.

## I. INTRODUCTION

High penetration of renewable generation, new energy storage devices, and emerging electricity market behavior have all thoroughly changed the characteristics of conventional power grids, bringing significant uncertainties (e.g., Californian duck curves) [1]. Conventional power grid operation methods that greatly rely on the system model now face grand challenges [2]. For example, if a mathematical model of the power grid is not available or is inaccurate, most of the existing model-based methods will lose their effectiveness. The same issue also occurs in optimal active power dispatch (OAPD) approaches where the system model is usually required [3].

When using conventional OAPD methods, the optimization problem is typically formulated as an objective function and is subjected to various constraints under known system conditions and models, resulting in an optimal solution derived from known system dynamics [4]. This type of formulation, however, has a major drawback. Once the real system condition changes and the system model becomes inaccurate, the calculated optimal operation cost may greatly deviate from the actual optimal cost. Furthermore, under certain circumstances, such as *N*-1 or *N-k*, the typical mathematical modeling approach will become even more challenging due to the high complexity of the problem [5].

To date, several adaptive methods have been proposed to handle a range of environment uncertainties [6-8]. In [6], an adaptive robust optimization method is proposed for multi-period economic dispatch (ED) to address the dynamic uncertainty caused by wind power. The highly dynamic uncertainty is formulated as a bounded variable, and the power flow constraint is simply formulated as a power balance constraint between load and generation. Reference [7] further supplements this method by presenting an adaptive robust ED method for tie-line scheduling while factoring in multi-area power systems with wind power penetration. A similar adaptive robust strategy for ED is proposed in [8], in which ED is implemented in a distributed manner using a consensus-based framework. Here, communication defects are studied during problem formulation. The adaptive-based methods are usually designed via a two-stage decision-making structure, with one of the stages estimating the time-varying system dynamics. However, a common issue existing in this approach is that the problems are formulated based on linearized system models derived from a dc power flow analysis – an assumption that is improper, particularly since the power system is well-known for its high nonlinearity and non-convexity, and especially when the system has high penetration of renewable generation and power electronic devices.

One solution proposed in [9] is to use a nonlinear function or neural network (NN) instead of an adaptive function to estimate the nonlinear system dynamics. However, the online-training process of this approach is computationally intensive and vulnerable to environmental noises. Even though it has been fully tested in the experimental environment, its practical performance cannot be guaranteed.

Recently, artificial intelligence (AI) in various applications such as AlphaGo, ATARI games, and robotics have shed light on potential use in deep reinforcement learning (DRL) approaches to solve optimal control problems [10]. After adequate off-line training, the DRL agent can make a series of optimal decisions to achieve the best goal.

The project is supported by State Grid Science and Technology Program.

In this paper, a DRL-based method is proposed to solve the AC-OAPD problem. We modify the double-deep-Q-network (DDQN) based algorithm proposed in [11] to dynamically adjust the generation output of each generator. The proposed method considers full ac power flow equations; similarly, the corresponding constraints are addressed in either the DRL formulation or the power system environment. It should be noted that since this work is mainly focused on the conceptual proof of the feasibility of using DRL to solve the OAPD problem, the thermal limit of the transmission line and voltage adjustments for the generator bus are not taken into consideration.

After a certain period of training, the DRL agent is able to solve the AC-OAPD problem quickly even under unseen scenarios. The proposed framework is a key component of the Grid Mind platform, which has been developed for autonomous grid operation using state-of-the-art AI techniques [10]. We believe that this framework holds promise to be generalized in the future to solve the AC-OPF problem. The effectiveness of the proposed DRL-based AC-OAPD method is demonstrated through simulations on the IEEE 14-bus system.

The remainder of this paper is organized as follows. Section II provides background on DRL and the formulation of AC-OAPD problem. Section III discusses in-depth implementation. In Section IV, case studies are performed using the proposed method, with promising results on the IEEE 14-bus test system. Finally, conclusions are drawn in Section V and future research work is also identified.

## II. METHODOLOGY OF DRL

### A. DRL: Background

A general framework diagram of DRL is illustrated in Fig. 1. The environment block is a physical or dynamic system, e.g., a power grid. The DRL agent is a controller that continuously interacts with the environment. The iteration of interactions starts when the DRL agent receives the state measurement ($s$) of the environment. The agent then provides a control action ($a$) according to the objective and definition of the reward function. After the environment conducts the action, it generates a new state ($s'$) and the corresponding reward ($r$). Based on the new observation of state and reward, the DRL agent then updates the framework and continues to improve its strategy.

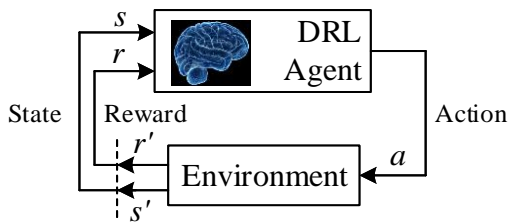

Fig. 1 The framework diagram of deep reinforcement learning.

Even though the global optimal solution of DRL cannot be guaranteed, the DRL mechanism makes it a promising tool to solve the optimization problem. For example, for a convex problem, the DRL solution can be guaranteed to be globally optimal. For a non-convex problem, on the other hand, various methods such as evolving program and random perturbation can be applied to avoid being trapped into the local optimal [16]. In this work, the double-deep-Q-network (DDQN) method with decaying ε-greedy policy is developed to solve the OAPD problem.

### B. Principles of DQN and DDQN

The deep-Q-network (DQN) is a value-based DRL algorithm developed from the classical RL method, i.e., Q-learning, where a deep neural network (DNN) is applied to estimate Q-function. Due to the estimation/prediction capability of DNN, the DQN method can handle the unseen states [12]. A general formulation for updating the value function can be presented as

$$Q'_{(s,a)} = Q_{(s,a)} + \alpha[r + \gamma max Q_{(s',a')} - Q_{(s,a)}] \quad (1)$$

where $\alpha$ is the learning rate and $\gamma$ is the discount rate; $s$ and $s'$ represent the current state and next state, respectively; $a$ and $a'$ represent the current action and next action, respectively. The parameters of DNN are updated by minimizing the error between actual and estimated Q-values ($[r + \gamma max Q_{(s',a')} - Q_{(s,a)}]$).

To solve the convergence problem, a double DNN structure is applied in DQN, which is then formulated as a new algorithm called DDQN. One of the DNNs, termed as the target network, fixes its parameters for a period of time and updates periodically. The other DNN called evaluation network continuously updates the parameters based on the estimation error of the target network. Both networks have the same structure, but with different parameters. In this way, the occasional oscillation caused by the big temporal difference error can be avoided.

To balance the exploration and exploitation, the decaying ε-greedy method is applied during the training process. Thus, in each iteration $i^{th}$, the DRL agent has a probability of $\varepsilon_i$ to choose an action randomly, and this probability keeps decaying, along with the training process as

$$\varepsilon_{i+1} = \begin{cases} r_d \times \varepsilon_i, \text{ if } \varepsilon_i > \varepsilon_{min} \\ \varepsilon_{min}, \quad \text{else} \end{cases} \quad (2)$$

where $r_d$ is a constant decay rate. By properly tuning the value of ε, the DRL agent will have a better chance to reach the global optimal, e.g., adding more perturbations in local optimal.

Various regularization methods have also been applied during the DRL agent training process to improve performance, e.g., random layer dropout, feature/batch normalization, and prioritized sampling [13].

## III. FORMULATION OF DRL-BASED OAPD

In this section, the DRL method is formulated to solve the optimal active power dispatch problem. First, the formulation of using proposed DRL to solve OAPD is presented. Then, the implementation details of the training environment of the DRL agent is illustrated.

### A. Formulation of DRL for OAPD

In the formulation, the full AC power flow equations are considered including the constraints of bus voltage, active power, and reactive power of generators. The bus voltage constraint is considered through training of DRL agents, which is similar to the autonomous voltage control function in the Grid Mind platform [10]. The active power generation limit of a generator is defined in the action space, and the reactive power generation limit of a generator is enforced by the in-house developed power flow solver for the Linux operating system. Since this work is mainly for conceptual proof, the line flow constraint is omitted due to the limitation of the current power flow solver in the current version; however, it can be easily incorporated into future versions.

The system states feeding into the DRL agent include the phasor voltage, the active, and reactive power flow on the transmission lines. However, for the generator located at the slack bus, which is used to balance the demand and supply, the action space of each generator is defined as [+0.5, -0.5, 0] MW. This means that the active power output of each generator can increase by 0.5 MW, decrease by 0.5 MW or remain unchanged at the current operating point. Thus, for a power system containing $N$ generators (excluding slacking bus), the total action space is $3^{N-1}$. Before the action is actually applied to the environment, the active power generation constraint of each generator is checked to guarantee that it is reasonable. In addition, the initial bus voltage magnitude settings of the generator are fixed at the optimal solution provided by the classic AC-OPF solver, i.e., MATPOWER [14]. Then, during the training process, the DRL agent continually checks to see if there is any voltage violation. If so, corresponding rewards are assigned to the agent to avoid future voltage violation. The predefined operation zone for voltage is illustrated in Fig. 2, where the desired target is defined to keep the voltage within 0.9 to 1.1 per unit (p.u.). Meanwhile, the agent will act to decrease the operation cost. Correspondingly, the reward function can be defined as

$$R_i = \begin{cases} (R_p - \beta * cost), \forall V_j \in [0.95,1.05] \text{pu} \\ (-R_n), \exists V_j \notin [0.95,1.05] \text{pu} \\ (-R_e), \text{Power flow diverged} \end{cases} \quad (3)$$

where $R_p$ is a positive reward, $-R_n$ is a negative reward, and $-R_e$ is a large negative penalty; $\beta$ is a constant weight to scale the generation cost.

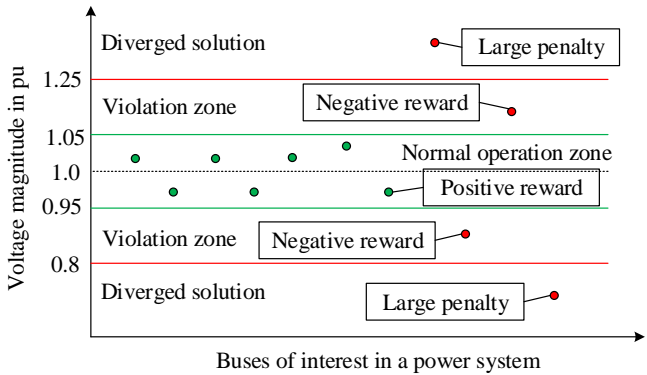

Fig. 2 The predefined voltage profile zone.

Define $P_i$ as the active power generation of generator $i$, a quadratic cost function that can be formulated as

$$cost = \sum_{i=1}^{n} a_i P_i^2 + b_i P_i \quad (4)$$

where $a_i$ and $b_i$ are cost coefficients. Note that since the system security is one of the primary concerns for the system operators, the three items ($R_p - \alpha$, $-R_n$ and $-R_e$) need to be properly weighted to achieve the desired target.

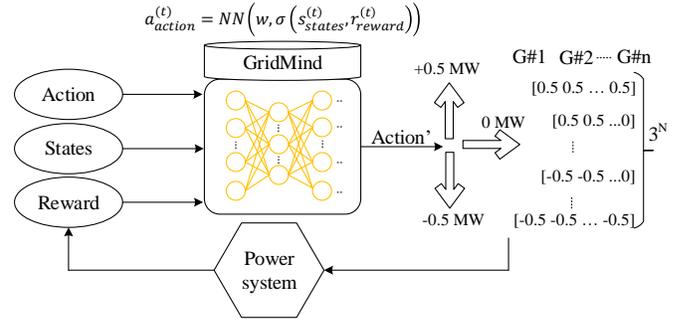

Fig. 3. Designed RL agent for optimal active power dispatch.

The diagram of implementing the proposed DRL agent for OPED is shown in Fig. 3. During the training process, a complete episode is defined as follows. The agent observes the state of an independent snapshot of power flow case and takes various actions to interact with the environment to maximize the reward until one of the termination conditions is met. The termination of one episode is triggered when one of the following conditions is met: 1) maximum iteration time is reached, 2) bus voltage violation occurs, and 3) power flow diverges. During the testing period, two different principles can be applied to obtain the optimal solution. The first method is to continuously compare the OAPD result with the benchmark system available in MATPOWER [14]. Once the calculated cost is equivalent to or fewer than the optimal cost, the DRL agent can output the solution. The second method is to allow the agent to take a predefined maximum iteration time. Once the maximum iteration time is reached, the result with the minimum cost will be the output.

In this paper, the second method is applied during the training period so that the DRL agent is encouraged to find a better solution, and the first method is applied during the testing period to evaluate the efficiency. In the next section, we show that although it takes hundreds of steps for the DRL agent to find the optimal solution during the training period, it only takes tens steps to solve the problem during the testing period.

*B. Implementation of DRL Agent for OAPD*

The platform used to train and test DRL agents for autonomous voltage control is Ubuntu 16 Linux Operation System (64 bit). This server is equipped with Intel Xeon E7-8893 v3 CPU at 3.2 GHz and 528 GB memory. All the DRL training and testing processes are performed on this platform.

To mimic a real power system environment, an in-house developed power grid simulator is adopted, which can be run in a Linux environment with multiple functional modules, such as AC power flow and voltage stability analysis. In this work, the AC power flow module is leveraged to interact with the DRL agent. Intermediate files are used to exchange information between grid simulator and the DRL agent, including power the flow information file saved in PSS/E v26 raw format and power flow solution results saved in text files. The proposed framework is programmed using Python language with TensorFlow libraries [15].

## IV. CASE STUDIES AND DISCUSSION

The proposed DRL method for autonomous voltage control is tested on the IEEE 14-bus system model. The single-line diagram is shown in Fig. 4. The IEEE 14-bus system consists of 14 buses, 5 generators, 11 loads, 17 lines, and 3 transformers. Since there is one generator located at the slack bus, the total action space is $3^4 = 81$. The total system load in the base case is 242 MW and 73.5 MVAr.

To demonstrate the effectiveness and adaptability of the proposed OAPD method, two different case studies are conducted. In the first study, the agent is trained using one base case on a normal operating condition. Then, the agent is tested on an additional 45 cases with the same loading condition, but different initial values of power generation. In the second case, the agent trained in the first case is tested on an additional 4 cases with different load levels, e.g., 80%~120%.

The detailed parameter settings of the tested system and the DRL agent are given in Table I.

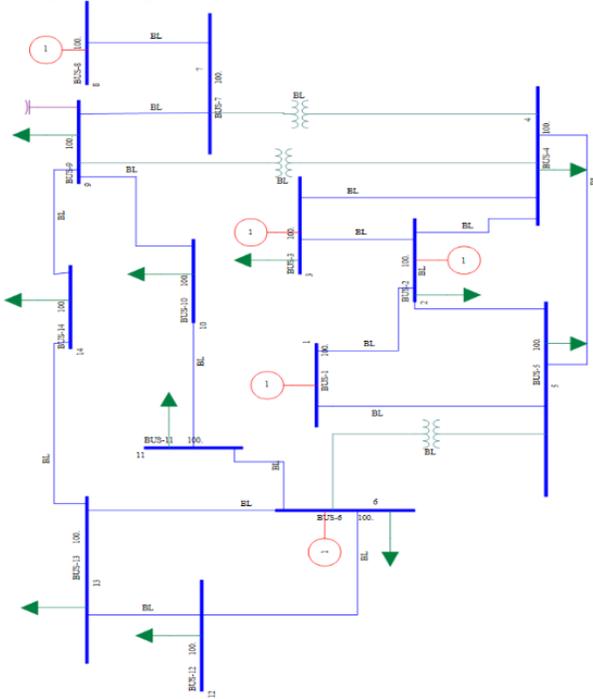

Fig. 4. Single-line diagram of the IEEE 14-bus system.

TABLE I
MAJOR PARAMETERS OF DRL

| Description | Parameters | Description | Parameters |
|---|---|---|---|
| $\alpha$ | 0.001 | $\gamma$ | 0.95 |
| $r_d$ | 0.999 | $\beta$ | 0.02 |
| Memory Size | 2000 | Mini Batch | 200 |
| P Constraint of G1 (MW) | 37.9~245.4 | P Constraint of G2 (MW) | 48.1~157.5 |
| P Constraint of G3 (MW) | 5.8~82.5 | P Constraint of G4 (MW) | 11.4~110.5 |
| P Constraint of G5 (MW) | 0.3~80.1 | Q Constraint of G1 (MVAR) | -132.2~76.1 |
| Q Constraint of G2 (MVAR) | -76.9~36.1 | Q Constraint of G3 (MVAR) | -0.8~48.5 |
| Q Constraint of G4 (MVAR) | -19.1~19.4 | Q Constraint of G5 (MVAR) | -8.0~19.2 |
| $a_1, a_2, a_3, a_4, a_5$ | 0.043, 0.25, 0.01, 0.01, 0.01 | $b_1, b_2, b_3, b_4, b_5$ | 20, 20, 40, 40, 40 |

### A. Case Study I – Different Initial Power Generation Points

In this case study, the loading condition is fixed. Under this load pattern, the optimal generation cost calculated with MATPOWER is 7.1341 k$/hr, which can be considered as a standard optimization solution. Representative results obtained by using DRL agent are summarized in Table II. As can be seen, based on the simulation results of 45 different cases, four different types of solutions are obtained. Among them, the standard optimization solution is achieved in 29 cases (i.e., 63%).

There are several reasons that lead to a sub-optimal solution. As mentioned before, the DRL algorithm is tented to reach the sub-optimal because of its gradient-based policy. Nevertheless, promising results are still observed by using different techniques. The second reason is the limitation of our in-house developed power flow solver. The minimum tolerance of power mismatch can only be set to 2e-3. With further development on the power flow solver, better performance can be achieved. The third reason is that the action can only be adjusted discretely in a step range of 0.5 MW. Thus, certain values of active power generation are omitted. With a smaller step range, e.g., 0.1 MW, the DRL agent will have a better chance to reach the global optimal. However, this becomes a trade-off between speed and accuracy since more iterations will be expected with a smaller step change.

An example in Fig. 5 presents the exploration process of the DRL agent to minimize the cost. As can be observed, during the exploration process, the DRL agent is trapped in the sub-optimal point several times. By including the corresponding perturbation techniques, the DRL agent is then forced to adjust its actions and policy for further exploration and to ultimately achieve its goal.

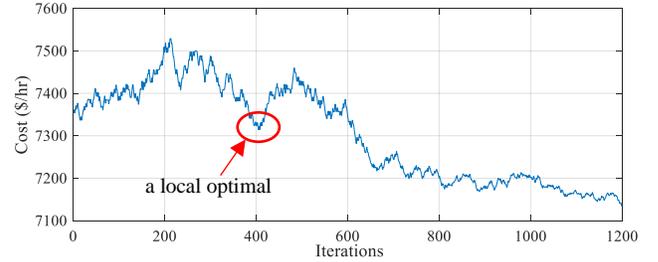

Fig. 5. Cost profile of DRL exploration process.

TABLE II
RESULT SUMMARY OF CASE STUDY I

| OAPD Solution of MATPOWER | | | |
|---|---|---|---|
| Case # | Active Power (MW) | Cost (k$/hr) | Power Mismatch |
| 1 | [176.9, 48.2, 5.8, 11.5, 0.3] | 7.1341 | 1e-9 |
| OAPD Solution of DRL | | | |
| Case # | Active Power (MW) | Cost (k$/hr) | Power Mismatch |
| 1 | [176.4, 48.2, 5.8, 11.9, 0.3] | 7.1361 | 4e-3 |
| 2 | [176.9, 48.2, 5.8, 11.5, 0.3] | 7.1341 | 2e-3 |
| 3 | [175.9, 48.2, 6.3 11.5, 0.8] | 7.1390 | 8e-3 |
| 4 | [176.2, 48.2, 5.8, 12.2, 0.3] | 7.1376 | 6e-3 |

It is also worth mentioning that during the training period, the DRL agent takes more than 1,000 steps to reach the optimal solution. However, during the testing period, the DRL agent only takes an average of 20 steps to find the solution. This further demonstrates that the strong learning and estimating capability of DRL mechanism make it an effective tool for solving the OAPD problem.

*B. Case Study II – Different Loading Conditions*

In this case, four different loading conditions of 80%, 90%, 110%, and 120%. are applied to test the effectiveness of the proposed method. Under each loading condition, the active power generation of each generator is first re-dispatched based on its inertia. Then, the DRL agent is applied to find the OAPD solution. The results are summarized in Table III.

As can be seen, even if the DRL agent is only trained in a base case, it can successfully solve the optimal generation dispatch problem under unseen scenarios. Moreover, the DRL agent is occasionally able to find a better solution, e.g., 120% loading condition (with larger power mismatch error), thereby further demonstrating the adaptability of the proposed DRL method for unknown system conditions.

TABLE III
RESULT SUMMARY OF CASE STUDY II

| OAPD Solution of MATPOWER | | | |
|---|---|---|---|
| Loading | Active Power (MW) | Cost (k$/hr) | Power Mismatch |
| 80% | [126.7, 48.1, 5.8, 11.4, 0.3] | 5.4663 | 1e-9 |
| 90% | [151.7, 48.1, 5.8, 11.4, 0.3] | 6.2656 | 1e-9 |
| 110% | [201.7, 48.2, 5.8, 11.5, 0.3] | 8.0338 | 1e-9 |
| 120% | [195.3, 48.1, 26.2, 11.4, 9.7] | 8.9876 | 1e-9 |
| OAPD Solution of DRL | | | |
| Case # | Active Power (MW) | Cost (k$/hr) | Power Mismatch |
| 80% | [125.5, 49.2, 5.8, 11.5, 0.3] | 5.4821 | 4e-3 |
| 90% | [151.6, 48.2, 5.8, 11.4, 0.3] | 6.2667 | 2e-3 |
| 110% | [201.7, 48.2, 5.8, 11.5, 0.3] | 8.0338 | 2e-3 |
| 120% | [215.4, 48.2, 8.3, 12.5, 6.3] | 8.9345 | 5e-2 |

In conclusion, all the above presented results have demonstrated the effectiveness of the proposed DRL-based OAPD method. Upon sufficient training, the DRL agent is shown to have the capability of solving the complex nonconvex optimization problem under unknown system dynamics.

## V. CONCLUSIONS AND FUTURE WORK

To effectively solve the optimal active power dispatch problem under growing uncertainties, this paper develops a novel method that uses the deep reinforcement learning method to search for optimal operating conditions with minimum generation cost. The detailed formulation and implementation process are introduced. The advantages of this approach, as well as its limitations are discussed and analyzed. Based on the results of tests performed on the IEEE 14-bus system, the proposed DRL-based method is demonstrated to be an effective tool to solve the OAPD problem under complex and unknown system conditions.

In future work, other practical constraints such as line flow limits will be considered and formulated into the DRL method. The voltage magnitude will be cooperatively adjusted, together with the reactive power. In addition, larger power system networks with or without contingencies will be used to test the performance of the DRL agent.